\documentclass[11pt,leqno]{amsart}
\usepackage{amssymb,amsmath,amsthm,latexsym}
\usepackage{graphicx}
\usepackage{epsfig}
\newtheorem{theorem}{Theorem}[section]
\newtheorem{lemma}[theorem]{Lemma}

\theoremstyle{remark}
\newtheorem{definition}[theorem]{Definition}
\newtheorem{remark}[theorem]{Remark}
\newtheorem{example}[theorem]{Example}

\newcommand\eL{{\mathcal L}}
\newcommand\M{{\mathcal M}}
\newcommand\R{{\mathcal R}}

\begin{document}

\title{On a homology of ternary groups with applications to knot theory}
\date{May 24, 2018}
\author{Maciej Niebrzydowski}
\address[Maciej Niebrzydowski]{Institute of Mathematics\\ 
Faculty of Mathematics, Physics and Informatics\\
University of Gda{\'n}sk, 80-308 Gda{\'n}sk, Poland}
\email{mniebrz@gmail.com}

\keywords{ternary group, homology, fundamental group, ternary knot group, Reidemeister moves, degenerate subcomplex, quasigroup, skew element, cocycle invariant}
\subjclass[2000]{Primary: 57M27, 20N10, 55N35; Secondary: 57Q45}

\thispagestyle{empty}

\begin{abstract}
We define a homology for ternary groups using both associativity and skew elements. We describe the odd-even construction which yields many examples of ternary groups. We define the ternary knot group, consider its homomorphisms into ternary groups, and discuss the applications.
\end{abstract}

\maketitle

\section{Introduction}

In this paper we describe a homology of ternary groups as an important example of ternary homology developed in \cite{Nie17a, Nie17b}. The presence of skew elements in ternary groups, together with repeated letters in the ternary differentials, leads to a major simplification and concise formulas. 
The ternary group of a knot that we are going to use is a modification of the classical knot group defined using Dehn presentation. 
We define it in order to have a structure from which to take homomorphisms into ternary groups, thus obtaining computable knot invariants. Ternary groups have a long history, see for example \cite{Post} for a thorough introduction to the subject of $n$-ary groups.
We begin with some preliminary definitions.

\begin{definition}\label{groupoid}
A {\it ternary groupoid} is a non-empty set $X$ equipped with a ternary operation $[\ ]\colon X^3\to X$. It is denoted by $(X,[\ ])$.
\end{definition}

\begin{definition}\label{semigroup}
We say that a ternary groupoid $(X, [\ ])$ is a {\it ternary semigroup} if the operation $[\ ]$ is associative, that is:
\[
[[abc]de]=[a[bcd]e]=[ab[cde]].
\]
\end{definition}

\begin{definition}\label{quasigroup0}
A {\it ternary quasigroup} is a ternary groupoid $(X, [\ ])$ such that for any quadruple $(a,b,c,d)$
of elements of $X$ satisfying $[abc]=d$, specification of any three elements of the quadruple determines the remaining one uniquely.
\end{definition}

For a wealth of material on $n$-ary quasigroups, see \cite{Bel} and \cite{BelSan}.

\begin{definition}\label{tgroup}
We say that a ternary groupoid $(X, [\ ])$ is a {\it ternary group} if it is both a ternary semigroup and a ternary quasigroup.
\end{definition}

\begin{definition}\label{skew}
In a ternary group $(X, [\ ])$, the unique solution $x$ of an equation $[aax]=a$ is called the {\it skew element} and is denoted by $\bar{a}$.
\end{definition}

\begin{lemma}[\cite{Dornte,DGG77}]
In any ternary group $(X, [\ ])$, for all $a$, $b$, $c\in X$, the following equalities hold:
\begin{align*}
&[ba\bar{a}]=[b\bar{a}a]=[a\bar{a}b]=[\bar{a}ab]=b,\\
&\overline{[abc]}=[\bar{c}\bar{b}\bar{a}],\\
&\bar{\bar{a}}=a.
\end{align*}
\end{lemma}

\begin{definition}
Let $(X, [\ ])$ be a ternary group. If the operation $[\ ]$ can be written in the form 
\[
[x_1x_2x_3]=x_1\circ x_2\circ x_3,
\]
for some binary group $(X,\circ)$, then we say that $(X, [\ ])$ is {\it derived} from $(X,\circ)$, or that $(X, [\ ])$ is {\it reducible} to $(X,\circ)$.
\end{definition}

In a ternary group derived from a binary group $(X,\circ)$, the skew element
$\bar{x}$ coincides with the inverse element $x^{-1}$ in $(X,\circ)$, see \cite{DG08}.

\section{The odd-even construction}

\begin{figure}
\begin{center}
\includegraphics[height=2.5 cm]{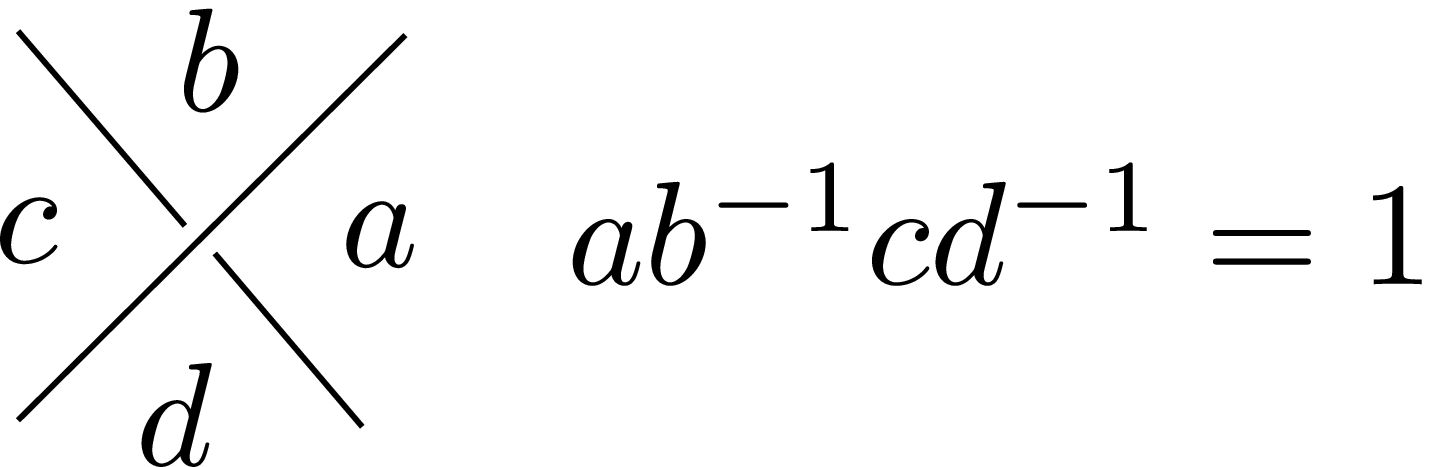}
\caption{}\label{Dehnop}
\end{center}
\end{figure}

The fundamental group of the complement of a knot in $\mathbb{R}^3$ can be given the following presentation, called Dehn presentation: generators are assigned to the regions in the complement of a knot diagram $D$ on a plane, and relations correspond to the crossings and are as in Figure \ref{Dehnop}. One of the generators, say the one corresponding to the unbounded region, is set equal to identity. Geometrically, a generator can be viewed as a loop originating from a fixed point $P$ above the diagram, piercing a region to which it is assigned, and returning to $P$ through a region labeled with the identity element. See, for example, \cite{Kauffor} for more details about this presentation. 

\begin{figure}
\begin{center}
\includegraphics[height=3.5 cm]{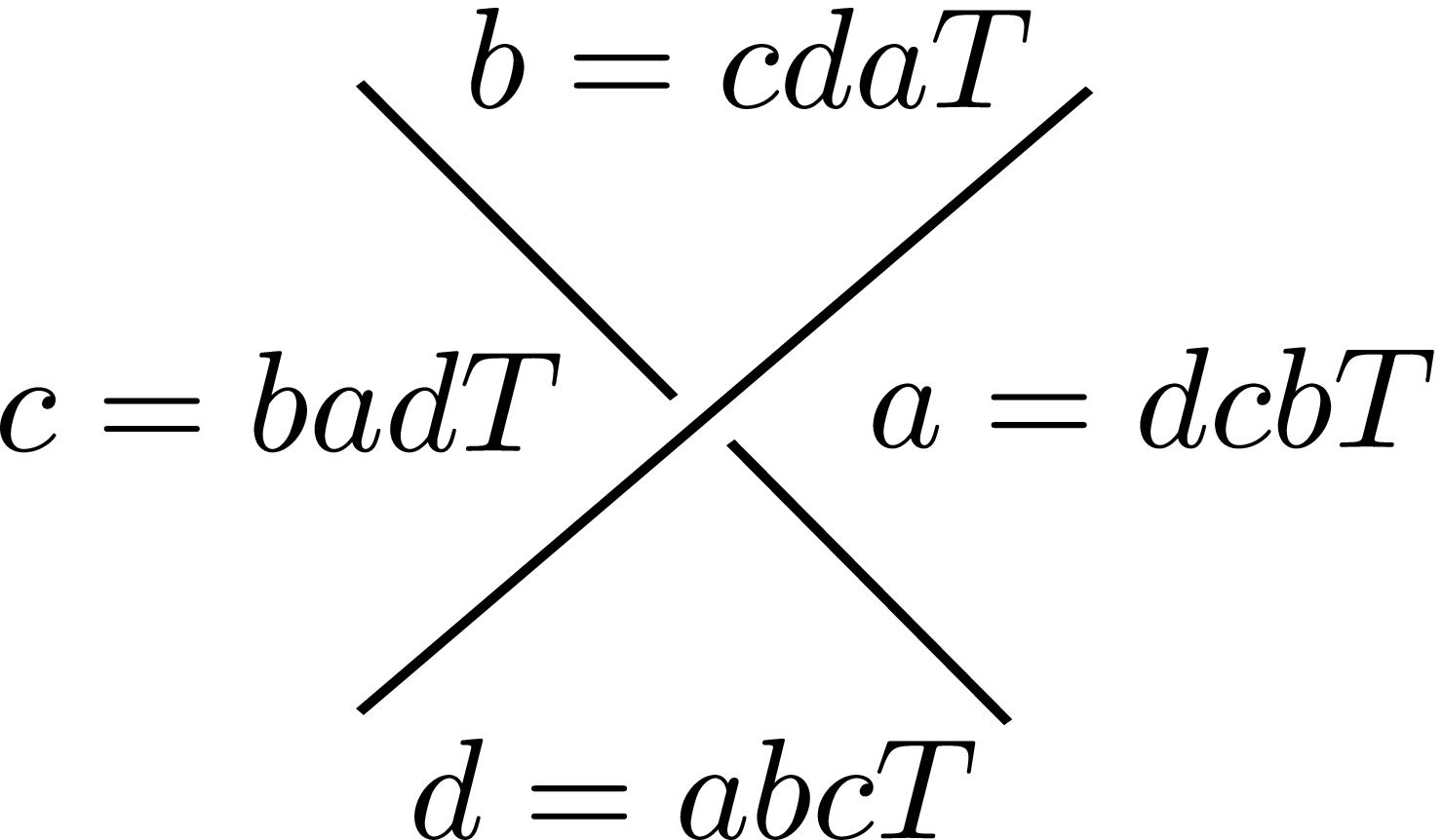}
\caption{}\label{ruleb}
\end{center}
\end{figure}

In \cite{Nie14}, we noted that the fundamental group relations between generators around a crossing can be realized using a ternary operation $xyzT=xy^{-1}z$. We see that each generator can be expressed using $T$ and the other three generators as in Fig. \ref{ruleb}. Namely, if $x$, $y$, $z$, $w\in X$ are the generators near a crossing, then $w=xyzT$, where $w$ and $x$ are assigned to the regions separated by an under-arc, and $x$, $y$ and $z$ are taken cyclically.

We will rewrite these relations on the level of ternary groups, but first we will describe a general method of finding ternary groups in some binary groups. The following is the motivating example for us. It provides examples of ternary groups not reducible to binary groups, although, as we can see from Theorem \ref{Posts}, every ternary group embeds in a binary group.

\begin{example}[\cite{Kasner}] \label{nonbin}
Let $S_n$ denote the symmetric group on $n$ elements, and let $A_n$ be the alternating group. There is a ternary group structure on $S_n\setminus A_n$,
defined by $[abc]=abc$. Here and later in the text, concatenation of letters without the presence of ternary brackets (as in $abc$) denotes the binary group multiplication. In particular, note that the binary identity element is removed.
\end{example}

More generally: let $G$ be a group with the set of generators $S$, and let $v$ be a map from the set of words in these generators (and their inverses) to $\mathbb{Z}_2$, sending words of even length to zero, and words of odd length to one. Suppose that this map is well defined on $G$, in the sense that if two words $w_1$ and $w_2$ represent the same group element $g\in G$, then
$v(w_1)=v(w_2)$. If $v(w_1)=0$, then we call such a $g$ {\it even}, if $v(w_1)=1$, we
say that $g$ is {\it odd} (with respect to the generators $S$). Then we can construct a ternary group $(X,[\ ])$ as in Example \ref{nonbin}: take $X$ to be the set of odd elements of $G$, and let $[abc]=abc$ be a triple product induced from $G$. Then, of course, odd elements are sent to odd elements by $[\ ]$, but not by the binary multiplication. We will call this method of obtaining ternary groups the {\it odd-even construction}. The set $H$ of even elements is also closed with respect to the ternary product $[abc]=abc$, but it contains the binary identity 1, is closed under the binary multiplication, and can be viewed as an ordinary binary group.

The above situation is reflected in the famous ``coset theorem" by Post about the existence of a binary group containing an $n$-ary group (here stated for ternary groups).

\begin{theorem}[\cite{Post}] \label{Posts}
For any ternary group $(X,[\ ])$, there exists a binary group $G$
and $H\lhd G$, such that $G/H\cong\mathbb{Z}_2$, and
\[
[xyz]=xyz,
\]
for all $x$, $y$, $z\in X$. Elements of $X$ can be identified with elements of $G$ that are not in $H$.
\end{theorem}  

In this paper we will use examples of ternary groups obtained by the odd-even construction via group presentations.

\begin{definition}
Let $G$ be an abstract group defined using presentation $P=\{S\ |\ R\}$, with the set of generators $S$ and the set of relations $R$. Suppose that the relations in $R$ do not change the parity of length of words, so that the parity map $v$ can be defined on group elements represented by words. To put it differently, if the relation $r\in R$ is written as $u=v$, for some words $u$ and $v$, then we require that the relator $uv^{-1}$ be of even length. Then, for odd elements $g_1$, $g_2$, $g_3\in G$ represented by words $w_1$, $w_2$, $w_3$ with an odd length, we have a ternary group product $[g_1g_2g_3]=w_1w_2w_3$, where $w_1w_2w_3$ is the concatenation of words. We will denote this ternary group of odd elements of $G$ (represented by words in the generators $S$ and their inverses) by $\mathcal{O}(P)$. We write $\mathcal{E}(P)$ for the often used binary group of even elements of $G$, with the operation $h_1h_2=v_1v_2$, for some words $v_1$ and $v_2$ with even length representing even elements $h_1$ and $h_2$. 
\end{definition}

\begin{figure}
\begin{center}
\includegraphics[height=3.5 cm]{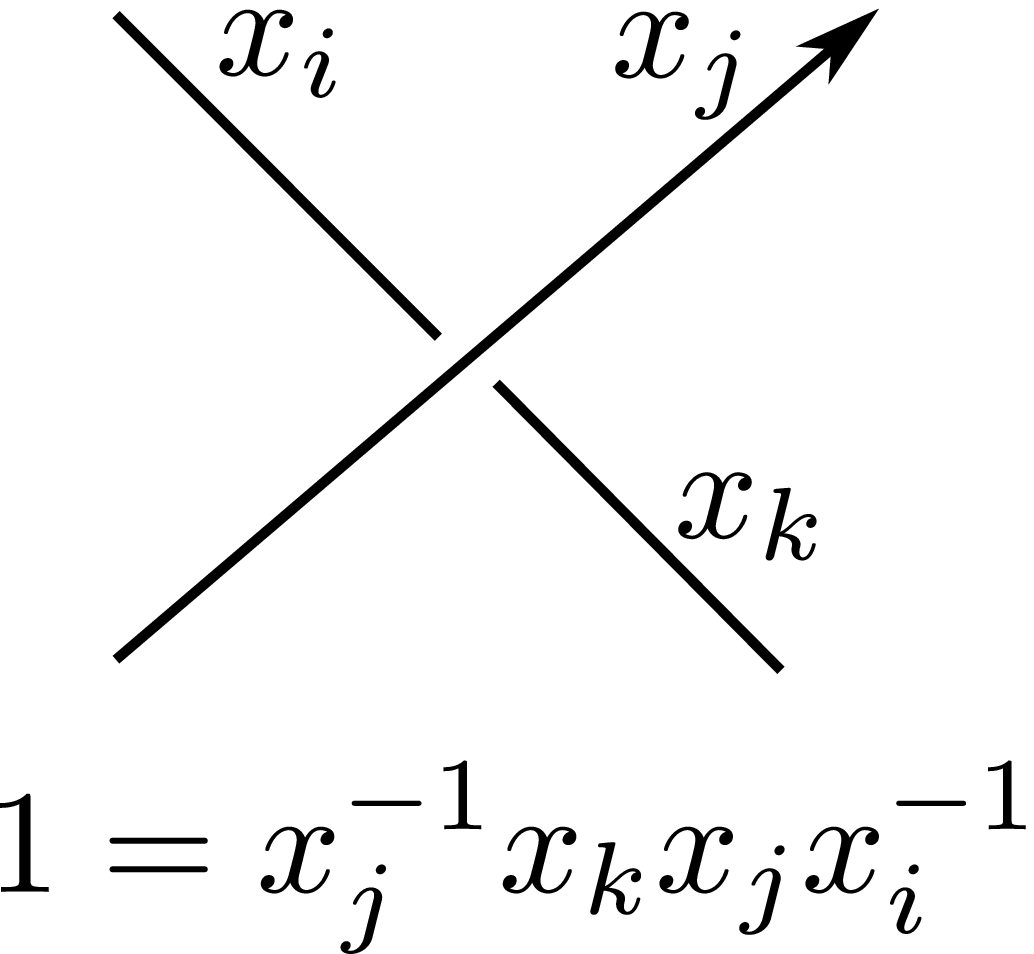}
\caption{}\label{wirt}
\end{center}
\end{figure}

\begin{example}
For a knot diagram $D$, let $\Pi_D$ denote the Wirtinger knot group presentation in which generators correspond to the arcs of $D$ and relations are assigned to crossings and are of the form: 
\begin{equation}\label{wirtrel}
1=x_j^{-1}x_kx_jx_i^{-1},
\end{equation}
as in Fig. \ref{wirt}. 
For $\Pi_D$ the parity map $v$ is well defined, because all the relators in this presentation have even length. 
Here $\mathcal{E}(\Pi_D)$ has an important topological interpretation: it is the fundamental group of the 2-fold cyclic covering of the complement of the knot in $\mathbb{S}^3$ (see \cite{Prz3col}). We are however interested in the ternary group $\mathcal{O}(\Pi_D)$. Note that the relations \eqref{wirtrel} can be rewritten in a ternary format: 
\begin{equation}\label{twirtrel}
x_i=[\bar{x}_jx_kx_j].
\end{equation}
The ternary group with presentation $T\Pi_D$ whose generators correspond to the arcs of $D$ and 
relations correspond to crossings and are of the form \eqref{twirtrel} is an invariant of the knot (up to isomorphism of ternary groups). 
\end{example}

\begin{example}
If we consider the braid group $B_n$ with the standard presentation
\begin{align*}
\{\sigma_1,\ldots,\sigma_{n-1}\ |\ &\sigma_i\sigma_j=\sigma_j\sigma_i\ \textrm{for}\ |i-j|\geq 2,\\
&\sigma_i\sigma_{i+1}\sigma_i=\sigma_{i+1}\sigma_i\sigma_{i+1}\ \textrm{for}\ i=1,\ldots,n-2\},
\end{align*}
then again the map $v$ is well defined on elements of $B_n$, because the relations do not change the parity of length of words. Therefore, we can take a triple product of braids with an odd length, thus obtaining a ternary group contained inside $B_n$.
We can define a ternary group $TB_n$ with presentation
\begin{align*}
\{\sigma_1,\ldots,\sigma_{n-1}\ |\ &[\sigma_i\sigma_{i+1}\sigma_i]=[\sigma_{i+1}\sigma_i\sigma_{i+1}]\ \textrm{for}\ i=1,\ldots,n-2,\\
&[\sigma_i\sigma_j\delta]=[\sigma_j\sigma_i\delta] \ \textrm{and}\ \\
&[\delta\sigma_i\sigma_j]=[\delta\sigma_j\sigma_i] \ \textrm{for}\ |i-j|\geq 2
\ \textrm{and an arbitrary generator}\ \delta\}.
\end{align*}
Here $\sigma_i$ corresponds to a crossing between the strings number $i$ and $i+1$, and $\bar{\sigma}_i$ corresponds to the opposite crossing.
Also note that equalities
\[
[\sigma_i\bar{\sigma}_i\delta]=[\bar{\sigma}_i\sigma_i\delta]=[\delta\sigma_i
\bar{\sigma}_i]=[\delta\bar{\sigma}_i\sigma_i]=\delta
\]
make sense geometrically, but do not need to be added to the set of relations, because they hold in any ternary group.
\end{example}

\begin{example}
More generally, the odd-even construction works for Artin group presentations
\begin{align*}
\{x_1,\ldots,x_{n}\ |\ &  <x_1,x_2>^{m_{1,2}}=<x_2,x_1>^{m_{2,1}},\ldots \\
&\ldots, <x_{n-1},x_n>^{m_{n-1,n}}=<x_n,x_{n-1}>^{m_{n,n-1}}\},  
\end{align*}
where $m_{i,j}=m_{j,i}\in\{2,3,\ldots,\infty\}$. Here $<x_i,x_j>$, for $m<\infty$, is the alternating product of $x_i$ and $x_j$ of length $m$, starting with $x_i$, for example:
\[ 
<x_i,x_j>^3=x_ix_jx_i.
\]
By convention, if $m=\infty$, then there is no relation for $x_i$ and $x_j$.
We can also define $\mathcal{O}(P)$ for a Coxeter group presentation
\begin{equation*}
P=\{x_1,\ldots,x_{n}\ |\ (x_ix_j)^{m_{ij}}=1 \},
\end{equation*}
where $m_{ii}=1$ and $m_{ij}\geq 2$ for $i\neq j$.
\end{example}

\begin{figure}
\begin{center}
\includegraphics[height=6 cm]{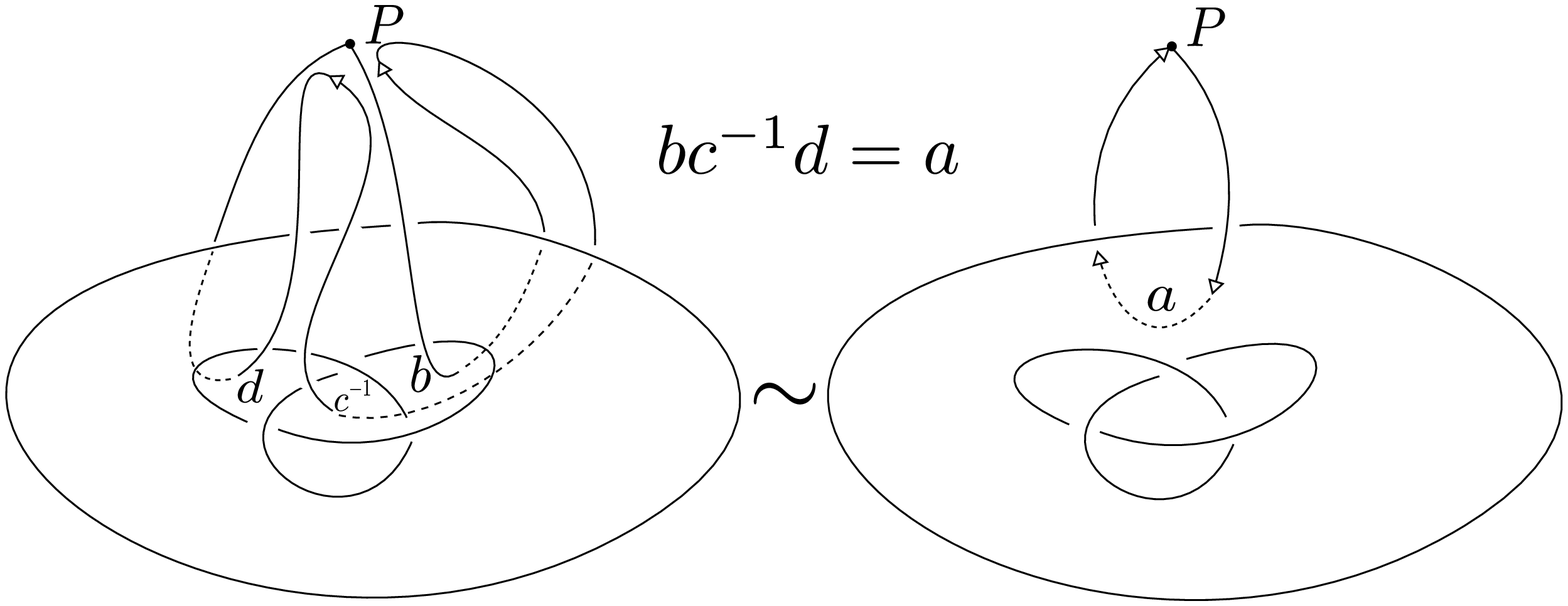}
\caption{}\label{relloops}
\end{center}
\end{figure}

\begin{figure}
\begin{center}
\includegraphics[height=3.5 cm]{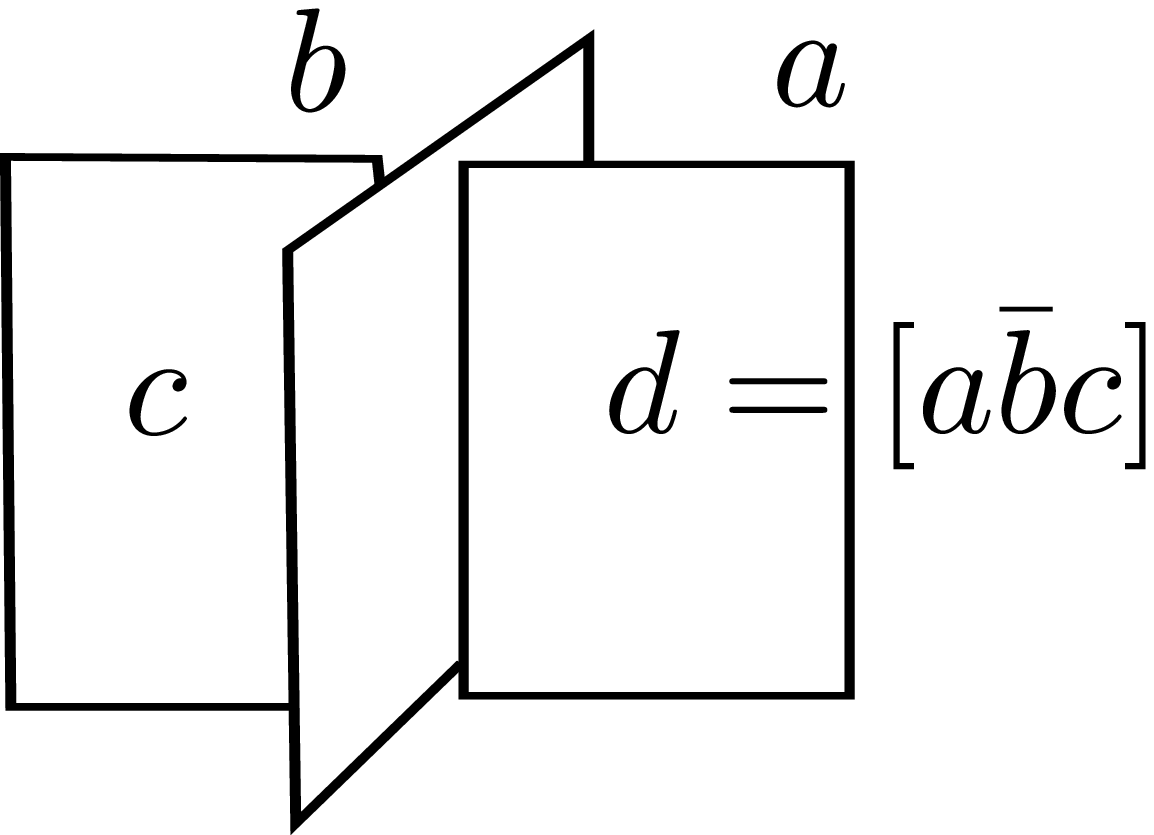}
\caption{}\label{surfdehn}
\end{center}
\end{figure}

\begin{example}
Finally, we consider the Dehn presentation of the knot group, assigned to a knot diagram $D$. Making one of the generators equal to identity adds a relation that spoils the parity map $v$, therefore we take the presentation without it, and denote it by $\pi_1(D)$. Then we can consider the ternary group $\mathcal{O}(\pi_1(D))$. 
We note that this minor change in the presentation could have a geometric meaning. Suppose that $D$ is on a disc $\mathbb{D}^2$, and a generator corresponds to a loop going from a point $P$ above $D$ through a given region and returning to $P$ underneath and outside the disc $\mathbb{D}^2$, thus linking itself with the boundary of the disc. With the ternary product being just the usual multiplication of three loops up to homotopy, the relations assigned to crossings are satisfied, and none of the generators corresponds to a trivial loop, see Fig. \ref{relloops}.

More abstractly, we define the {\it ternary knot group} as the ternary group with presentation $\pi_1^T(D)$ whose generators correspond to regions in the complement of $D$, and relations are assigned to crossings and are of the form $d=[a\bar{b}c]$ for a crossing as in Fig. \ref{Dehnop}. The operation $abcT=[a\bar{b}c]$ satisfies the axioms \eqref{first} and \eqref{second} (see \cite{NPZ18}), and it follows that the ternary knot group is a knot invariant (see \cite{Nie14, Nie17a, Nie17b} for more general theory of invariants obtained from ternary algebras satisfying \eqref{first} and \eqref{second}; in case of virtual knots, see \cite{NelPic18}).
Analogous ternary group can be assigned to a broken diagram of a knotted surface (see \cite{KnSurf} for the theory of such diagrams). In this case generators correspond to the regions in the complement of the projection of the surface in $\mathbb{R}^3$, and relations are assigned to the curves of double points and are as in Fig. \ref{surfdehn}. Then such a ternary group, up to isomorphism, is an invariant of Roseman moves. 
\end{example}

\begin{remark}
Let $G$ be a group with an odd-even construction. Then we can consider various notions of action by pairs of even elements on the set of odd elements, for example:
$(b_1,b_2,a)\mapsto b_1ab_2$, where $b_1$, $b_2$ are even, and $a$ is odd.
\end{remark} 

\section{Homology of ternary groups}

\begin{figure}
\begin{center}
\includegraphics[height=6 cm]{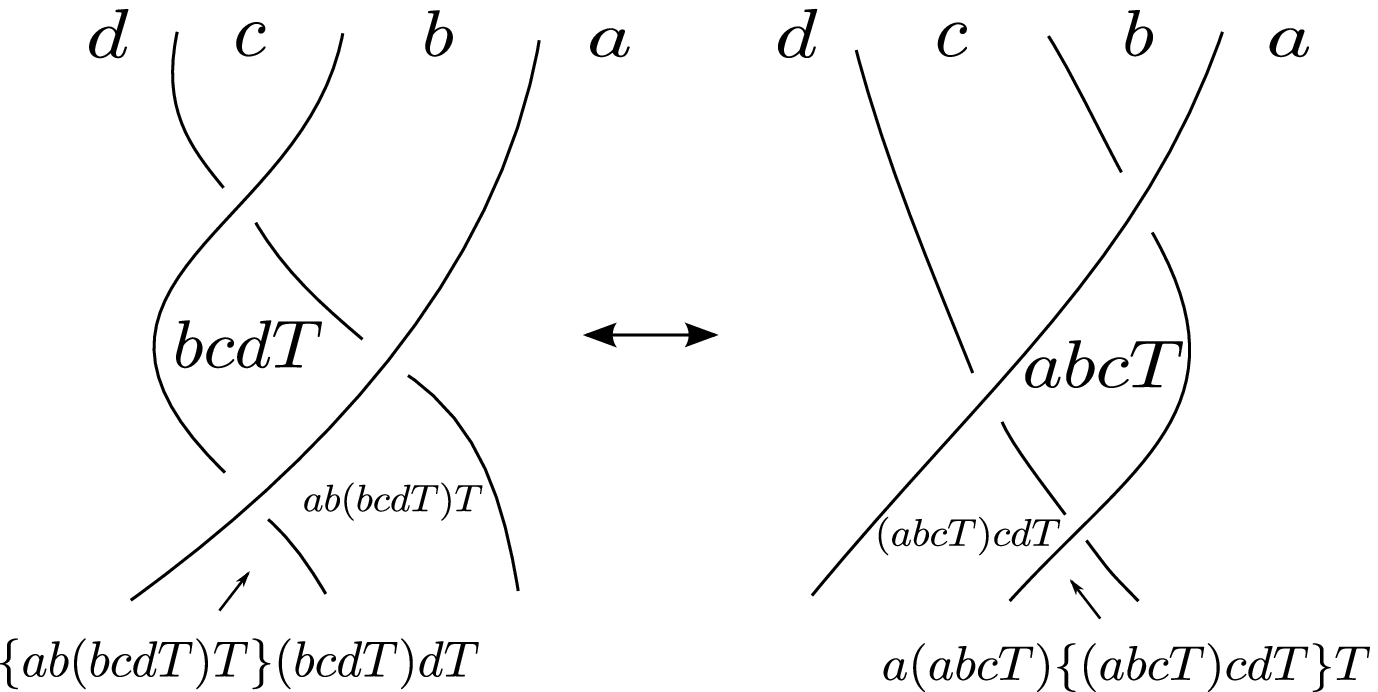}
\caption{}\label{R3}
\end{center}
\end{figure}

In \cite{Nie17a}, we defined homology for ternary algebras $(X,T)$ satisfying two axioms derived from colorings of regions and the third Reidemeister move, see Fig. \ref{R3}. In the Figure, we reversed the direction of taking arguments compared to our previous papers, to match exactly the conventions used in the literature for the Dehn presentation. The axioms are: 
\begin{equation} \label{first}
\forall_{a,b,c,d\in X} \quad (abcT)cdT=\{ab(bcdT)T\}(bcdT)dT
\end{equation}
and
\begin{equation} \label{second}
\forall_{a,b,c,d\in X} \quad ab(bcdT)T=a(abcT)\{(abcT)cdT\}T. 
\end{equation}
Here we used the notation $abcT$ instead of the bracket notation.

For a ternary group $(X, [\ ])$, and arbitrary elements $a$, $b$, $c\in X$, let 
\[
abcT=[a\bar{b}c].
\]
As we have already mentioned, such $T$ satisfies the axioms (\ref{first}) and (\ref{second}).
Now we define the chain groups. Let 
$C_n(X)=Z\langle X^{n+2}\rangle$ be the free abelian group generated by $(n+2)$-tuples $(x_0,x_1,\ldots, x_n,x_{n+1})$ of elements of $X$, for $n\geq -1$, and let $C_{-2}(X)=\mathbb{Z}$. We will not rewrite here all the differentials defined in \cite{Nie17a}, but we will indicate how the simplification for the operation
$abcT=[a\bar{b}c]$ is occuring. For example:
\begin{align*}
\partial_3(a,b,c,d,e)& =(b,c,d,e)-(a,abcT,(abcT)cdT,\{(abcT)cdT\}deT)\\
& -(abcT,c,d,e)+(a,b,bcdT,(bcdT)deT)\\
& +(ab(bcdT)T,bcdT,d,e)-(a,b,c,cdeT)\\
& -(ab\{bc(cdeT)T\}T,bc(cdeT)T,cdeT,e)+(a,b,c,d)
\end{align*}
becomes
\begin{align*}
\partial_3(a,b,c,d,e)& =(b,c,d,e)-(a,[a\bar{b}c],[[a\bar{b}c]\bar{c}d],[[[a\bar{b}c]\bar{c}d]\bar{d}e])\\
& -([a\bar{b}c],c,d,e)+(a,b,[b\bar{c}d],[[b\bar{c}d]\bar{d}e])\\
& +([a\bar{b}[b\bar{c}d]],[b\bar{c}d],d,e)-(a,b,c,[c\bar{d}e])\\
& -([a\bar{b}[b\bar{c}[c\bar{d}e]]],[b\bar{c}[c\bar{d}e]],[c\bar{d}e],e)+(a,b,c,d),
\end{align*}
and after simplifications using associativity and the properties of skew elements,
\begin{align*}
\partial_3(a,b,c,d,e)& =(b,c,d,e)-(a,[a\bar{b}c],[a\bar{b}d],[a\bar{b}e])\\
& -([a\bar{b}c],c,d,e)+(a,b,[b\bar{c}d],[b\bar{c}e])\\
& +([a\bar{c}d],[b\bar{c}d],d,e)-(a,b,c,[c\bar{d}e])\\
& -([a\bar{d}e],[b\bar{d}e],[c\bar{d}e],e)+(a,b,c,d).
\end{align*}
Now we write the general differential $\partial_n\colon C_n(X) \to C_{n-1}(X)$.
\begin{align*}
&\partial_n(x_0,x_1,\ldots,x_n,x_{n+1})=(x_1,\ldots,x_n,x_{n+1})\\
&+\sum_{i=1}^n(-1)^i\{([x_0\bar{x}_ix_{i+1}],[x_1\bar{x}_ix_{i+1}],\ldots,
[x_{i-1}\bar{x}_ix_{i+1}],\hat{x}_i,x_{i+1},\ldots,x_n,x_{n+1})\\
&+(x_0,x_1,\ldots,x_{i-1},\hat{x}_i,[x_{i-1}\bar{x}_ix_{i+1}],\ldots,
[x_{i-1}\bar{x}_ix_{n}],[x_{i-1}\bar{x}_ix_{n+1}])\}\\
&+(-1)^{n+1}(x_0,x_1,\ldots,x_n),
\end{align*}
where $\hat{x}_i$ denotes a missing element. We also take
\begin{align*}
& \partial_{-1}(x_0)=0,\\
& \partial_0(x_0,x_1)=x_1-x_0.
\end{align*}
Since the homology defined in \cite{Nie17b} is really a difference of two differentials, here too we can separate $\partial$ into two parts:
$\partial_n=\partial_n^L-\partial_n^R$, where
\begin{align*}
&\partial^L_n(x_0,x_1,\ldots,x_n,x_{n+1})=(x_1,\ldots,x_n,x_{n+1})\\
&+\sum_{i=1}^n(-1)^i([x_0\bar{x}_ix_{i+1}],[x_1\bar{x}_ix_{i+1}],\ldots,
[x_{i-1}\bar{x}_ix_{i+1}],\hat{x}_i,x_{i+1},\ldots,x_n,x_{n+1}),
\end{align*}
and
\begin{align*}
&\partial^R_n(x_0,x_1,\ldots,x_n,x_{n+1})=(-1)^{n}(x_0,x_1,\ldots,x_n)\\
&+\sum_{i=1}^n(-1)^{i+1} (x_0,x_1,\ldots,x_{i-1},\hat{x}_i,[x_{i-1}\bar{x}_ix_{i+1}],\ldots,[x_{i-1}\bar{x}_ix_{n}],[x_{i-1}\bar{x}_ix_{n+1}]).
\end{align*}

In connection with the third Reidemeister move, the differential of interest is:
\begin{align*}
\partial_2(a,b,c,d)=&(b,c,d)-([a\bar{b}c],c,d)-(a,[a\bar{b}c],[a\bar{b}d])\\
+&([a\bar{c}d],[b\bar{c}d],d)+(a,b,[b\bar{c}d])-(a,b,c).
\end{align*}

With the operation $abcT=[a\bar{b}c]$, $X$ is a quasigroup, so we can consider the left, middle and right divisions, denoted by $\eL$, $\M$, $\R$ as in \cite{Nie17b}. For our operation $T$, we find that
\begin{align*}
abc\eL&=[a\bar{c}b],\\
abc\M&=[c\bar{b}a],\\
abc\R&=[b\bar{a}c].
\end{align*}
In \cite{Nie17b} we defined a degenerate subcomplex for homology of ternary quasigroups satisfying \eqref{first} and \eqref{second}, corresponding to the first Reidemeister move. Here we recall that in this subcomplex 
$C_n^D(X)$ denotes the free abelian group generated by $(n+2)$-tuples 
$x=(x_0,x_1,\ldots, x_n,x_{n+1})$ of elements of $X$ satisfying at least one of the conditions:
\begin{itemize}
\item[(D1)] $x$ contains $a$, $b$, $abb\R$ on three consecutive coordinates, for some $a$ and $b\in X$;
\item[(D2)] $x$ contains $bba\eL$, $b$, $a$ on three consecutive coordinates, for some $a$ and $b\in X$. 
\end{itemize}
For $n<1$, we take $C_n^D(X)=0$. 

For the operation $abcT=[a\bar{b}c]$, this translates into:
$C_n^D(X)$ is the free abelian group generated by $(n+2)$-tuples 
$x=(x_0,x_1,\ldots, x_n,x_{n+1})$ of elements of $X$ containing, on three consecutive coordinates, a triple $a$, $b$, $[b\bar{a}b]$ or a triple
$[b\bar{a}b]$, $b$, $a$, for some $a$ and $b\in X$.
Then we have
\begin{theorem}[\cite{Nie17b}] 
\[\partial_n(C_n^D)\subset C_{n-1}^D.\]
\end{theorem}

\begin{definition}
We consider the quotient complex
\begin{equation*} 
(C_n^{N}(X),\partial_n)=(C_n(X)/C_n^D(X),\partial_n),
\end{equation*}
with induced differential (and the same notation). We define the {\it homology of a ternary group} $(X, [\ ])$ as the homology of this quotient complex, and denote it by $H_{*}^N(X)$.
\end{definition}

Of course, the above definition comes more from knot theory, than from algebraic considerations.

Using GAP \cite{GAP4}, we performed some calculations for presentations of triangle groups  
\[
\triangle(l,m,n)=\{ a,b,c\ |\ a^2=b^2=c^2=(ab)^l=(bc)^n=(ca)^m=1\},
\]
where $l$, $m$ and $n$ are integers greater than or equal to 2.
Note that in this case $\mathcal{E}(\triangle(l,m,n))$ is known as von Dyck group. Our calculations detected torsion in $H_1^N(\mathcal{O}(\triangle(l,m,n)))$ for some triples $(l,m,n)$ in the spherical case (i.e.
$1/l+1/m+1/n>1$) for which $\mathcal{O}(\triangle(l,m,n))$ is finite.

We begin with an explicit example of such ternary group. Consider $\triangle(2,2,3)$. We number its odd elements as follows:
1) $c$, 2) $a$, 3) $b$, 4) $a*c*b$, 5) $b*c*b$, 6) $a*b*c$. With this numbering, the multiplication cube for $\mathcal{O}(\triangle(2,2,3))$ can be sliced into the following six matrices, each matrix for a fixed first coordinate of the triple product $[x,y,z]=xyz$.
For example $[1,2,3]=4$ and $[2,3,4]=5$.

\[ 
\begin{array}{|c| c cccccc} 
1yz&1 & 2 & 3 & 4 & 5 & 6 \\
\hline 
1 & 1 & 2 & 3 & 4 & 5 & 6 \\
2 & 2 & 1 & 4 & 3 & 6 & 5 \\
3 & 5 & 4 & 1 & 6 & 3 & 2 \\ 
4 & 6 & 3 & 2 & 5 & 4 & 1 \\
5 & 3 & 6 & 5 & 2 & 1 & 4 \\ 
6 & 4 & 5 & 6 & 1 & 2 & 3 \\
\end{array}
\begin{array}{|c| c cccccc} 
2yz&1 & 2 & 3 & 4 & 5 & 6 \\
\hline 
1 & 2 & 1 & 4 & 3 & 6 & 5 \\
2 & 1 & 2 & 3 & 4 & 5 & 6 \\
3 & 6 & 3 & 2 & 5 & 4 & 1 \\
4 & 5 & 4 & 1 & 6 & 3 & 2 \\ 
5 & 4 & 5 & 6 & 1 & 2 & 3 \\ 
6 & 3 & 6 & 5 & 2 & 1 & 4 \\
\end{array}
\begin{array}{|c| c cccccc} 
3yz&1 & 2 & 3 & 4 & 5 & 6 \\
\hline 
1 & 3 & 6 & 5 & 2 & 1 & 4 \\
2 & 6 & 3 & 2 & 5 & 4 & 1 \\ 
3 & 1 & 2 & 3 & 4 & 5 & 6 \\ 
4 & 4 & 5 & 6 & 1 & 2 & 3 \\ 
5 & 5 & 4 & 1 & 6 & 3 & 2 \\ 
6 & 2 & 1 & 4 & 3 & 6 & 5 \\ 
\end{array}
\]

\[
\begin{array}{|c| c cccccc} 
4yz&1 & 2 & 3 & 4 & 5 & 6 \\
\hline 
1 & 4 & 5 & 6 & 1 & 2 & 3 \\
2 & 5 & 4 & 1 & 6 & 3 & 2 \\ 
3 & 2 & 1 & 4 & 3 & 6 & 5 \\ 
4 & 3 & 6 & 5 & 2 & 1 & 4 \\ 
5 & 6 & 3 & 2 & 5 & 4 & 1 \\ 
6 & 1 & 2 & 3 & 4 & 5 & 6 \\ 
\end{array}
\begin{array}{|c| c cccccc} 
5yz&1 & 2 & 3 & 4 & 5 & 6 \\
\hline 
1 & 5 & 4 & 1 & 6 & 3 & 2 \\ 
2 & 4 & 5 & 6 & 1 & 2 & 3 \\ 
3 & 3 & 6 & 5 & 2 & 1 & 4 \\ 
4 & 2 & 1 & 4 & 3 & 6 & 5 \\ 
5 & 1 & 2 & 3 & 4 & 5 & 6 \\ 
6 & 6 & 3 & 2 & 5 & 4 & 1 \\ 
\end{array}
\begin{array}{|c| c cccccc} 
6yz&1 & 2 & 3 & 4 & 5 & 6 \\
\hline 
1 & 6 & 3 & 2 & 5 & 4 & 1 \\
2 & 3 & 6 & 5 & 2 & 1 & 4 \\ 
3 & 4 & 5 & 6 & 1 & 2 & 3 \\ 
4 & 1 & 2 & 3 & 4 & 5 & 6 \\ 
5 & 2 & 1 & 4 & 3 & 6 & 5 \\ 
6 & 5 & 4 & 1 & 6 & 3 & 2 \\
\end{array}
\]
Note that the skew element operation is as follows: $\bar{1}=1$, $\bar{2}=2$, 
$\bar{3}=3$, $\bar{4}=6$, $\bar{5}=5$, $\bar{6}=4$.
Computer calculations show that for this ternary group
\[
tor H_1^N(\mathcal{O}(\triangle(2,2,3)))=\mathbb{Z}_9.
\]
Here are some other GAP results:
\begin{align*}
&tor H_1^N(\mathcal{O}(\triangle(2,2,2)))=0,\\
&tor H_1^N(\mathcal{O}(\triangle(2,2,4)))=\mathbb{Z}_2 + \mathbb{Z}_4,\\
&tor H_1^N(\mathcal{O}(\triangle(2,2,5)))=\mathbb{Z}_5 + \mathbb{Z}_{25},\\
&tor H_1^N(\mathcal{O}(\triangle(2,2,6)))=\mathbb{Z}_3 + \mathbb{Z}_9,\\
&tor H_1^N(\mathcal{O}(\triangle(2,3,3)))=\mathbb{Z}_4.
\end{align*}
Note that $\mathcal{O}(\triangle(2,2,2))$ has four elements, 
$\mathcal{O}(\triangle(2,2,4))$ has eight, $\mathcal{O}(\triangle(2,2,5))$ has ten, and both $\mathcal{O}(\triangle(2,2,6))$ and $\mathcal{O}(\triangle(2,3,3))$
have twelve elements.

\section{Homological invariants of knots}

Now we use homomorphisms from ternary knot groups to ternary groups as a basis for homological invariants of knots. The number of such homomorphisms to a given ternary group $(X, [\ ])$ is a knot invariant, but it is just the first step on the way to defining interesting knot invariants. For a comparison to and survey of invariants obtained from rack/quandle (co)homology, see \cite{Surf4}. 

\begin{figure}
\begin{center}
\includegraphics[height=2.5 cm]{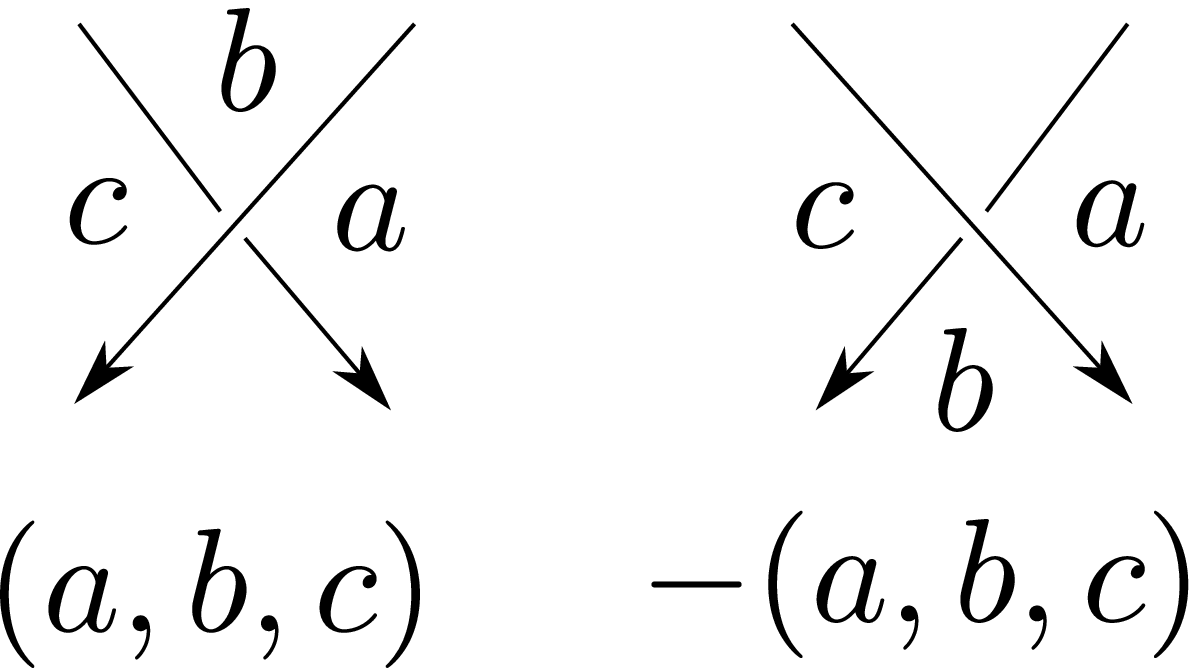}
\caption{}\label{chaindehn}
\end{center}
\end{figure}

\begin{figure}
\begin{center}
\includegraphics[height=4 cm]{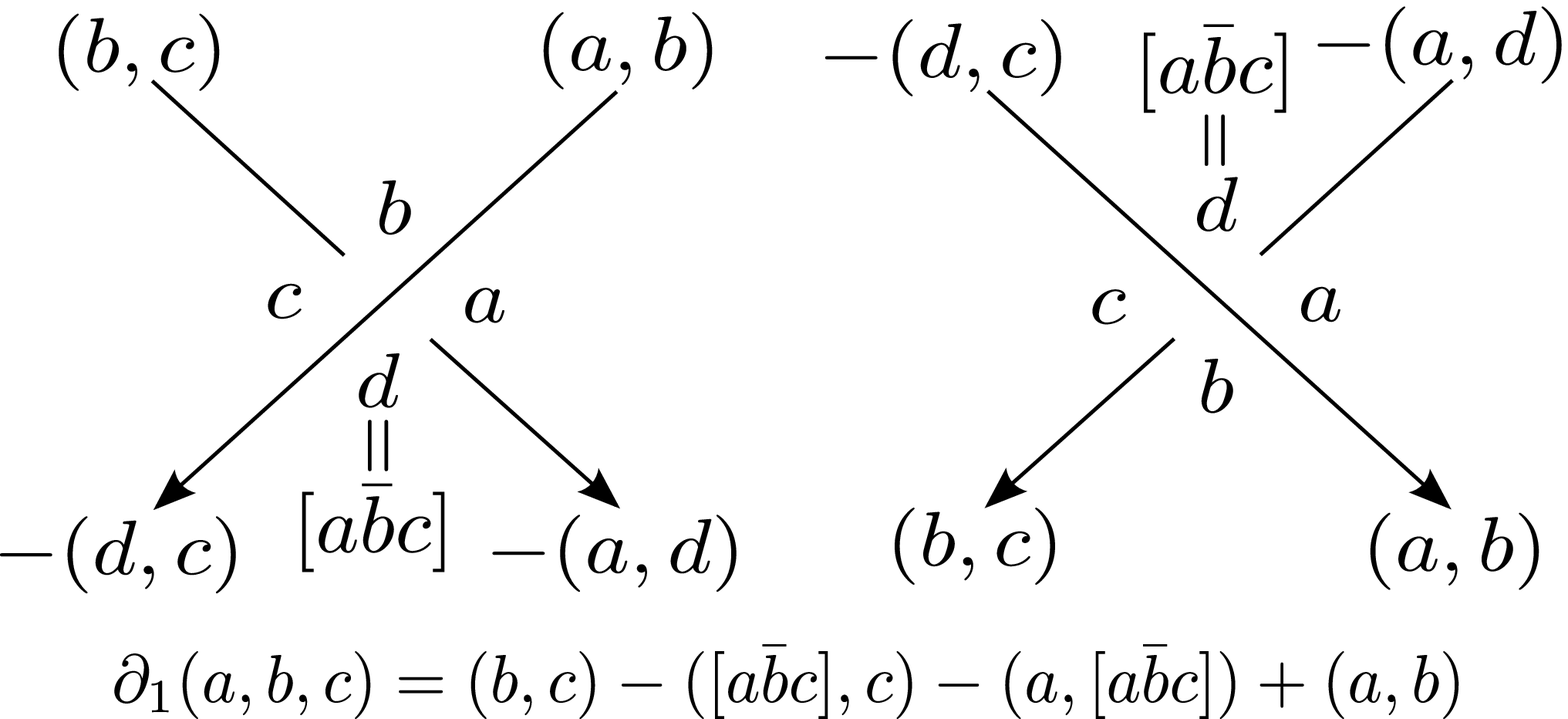}
\caption{}\label{cycleDehn}
\end{center}
\end{figure}

Let $D$ be an oriented knot diagram, and let $\gamma\colon \pi^T(D)\to X$ denote a homomorphism from a ternary knot group to a ternary group $(X, [\ ])$.
It gives a coloring of regions in the complement of the universe of $D$ by assigning a color $\gamma(x_i)$ to a region corresponding to a generator $x_i$.
To such colored $D$, we can assign a cycle in 
$H_1^N(X)$ by taking the sum of chains, as in Fig. \ref{chaindehn}, over all crossings of $D$. Fig. \ref{cycleDehn} shows that it really is a cycle: for an incoming edge of a positive crossing, the pair of colors surrounding it appears in $\partial_1(a,b,c)$ with a positive sign, the pairs corresponding to outgoing edges have a negative sign, so there is a reduction with terms at the next crossing where the given outgoing edge becomes incoming; the sign convention ensures that it also works for a negative crossing. We denote such a cycle by $c_\gamma(D)$.

\begin{figure}
\begin{center}
\includegraphics[height=5.5 cm]{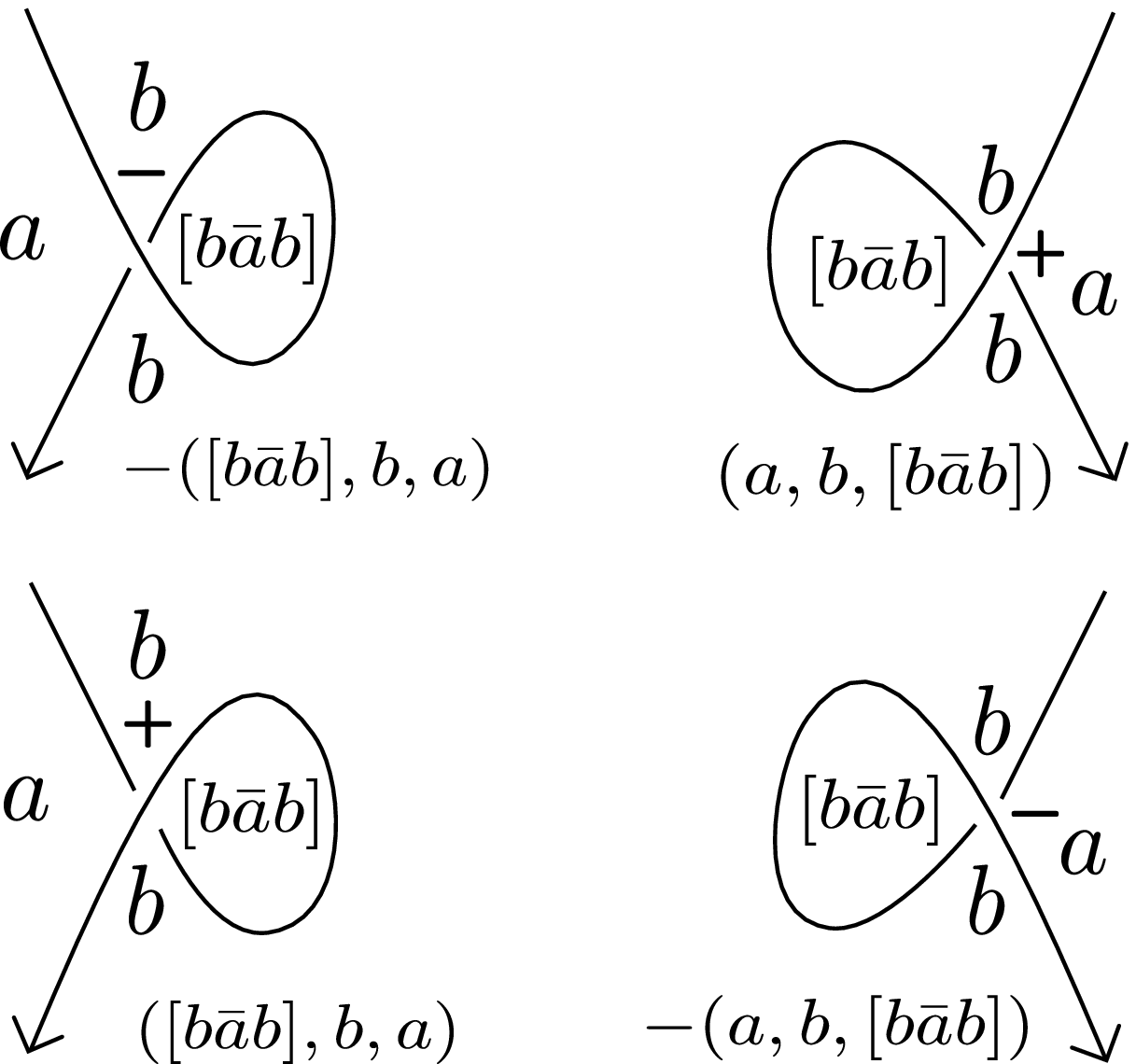}
\caption{}\label{R1dehn}
\end{center}
\end{figure}

\begin{figure}
\begin{center}
\includegraphics[height=4 cm]{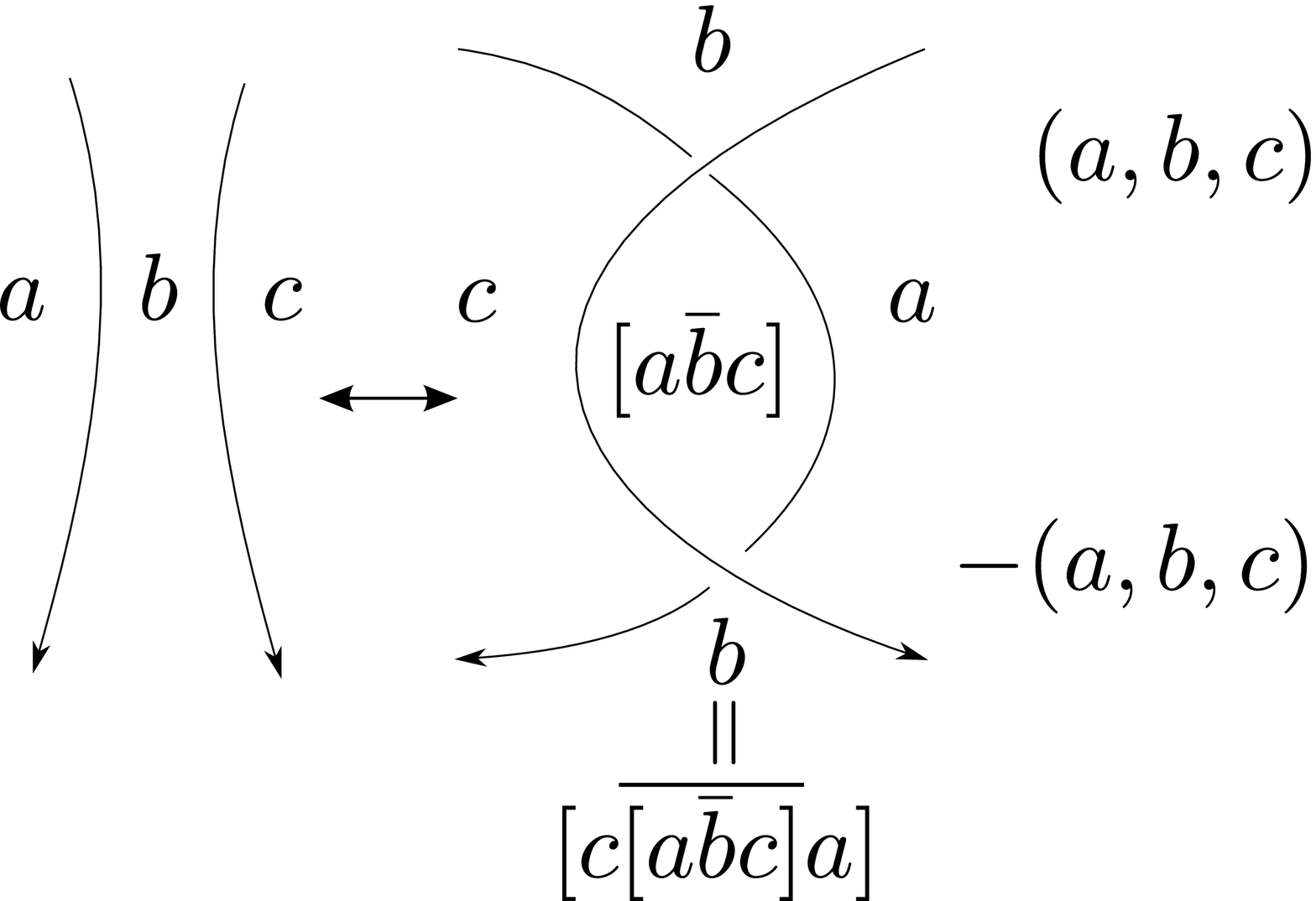}
\caption{}\label{seconddehn}
\end{center}
\end{figure}

\begin{figure}
\begin{center}
\includegraphics[height=5 cm]{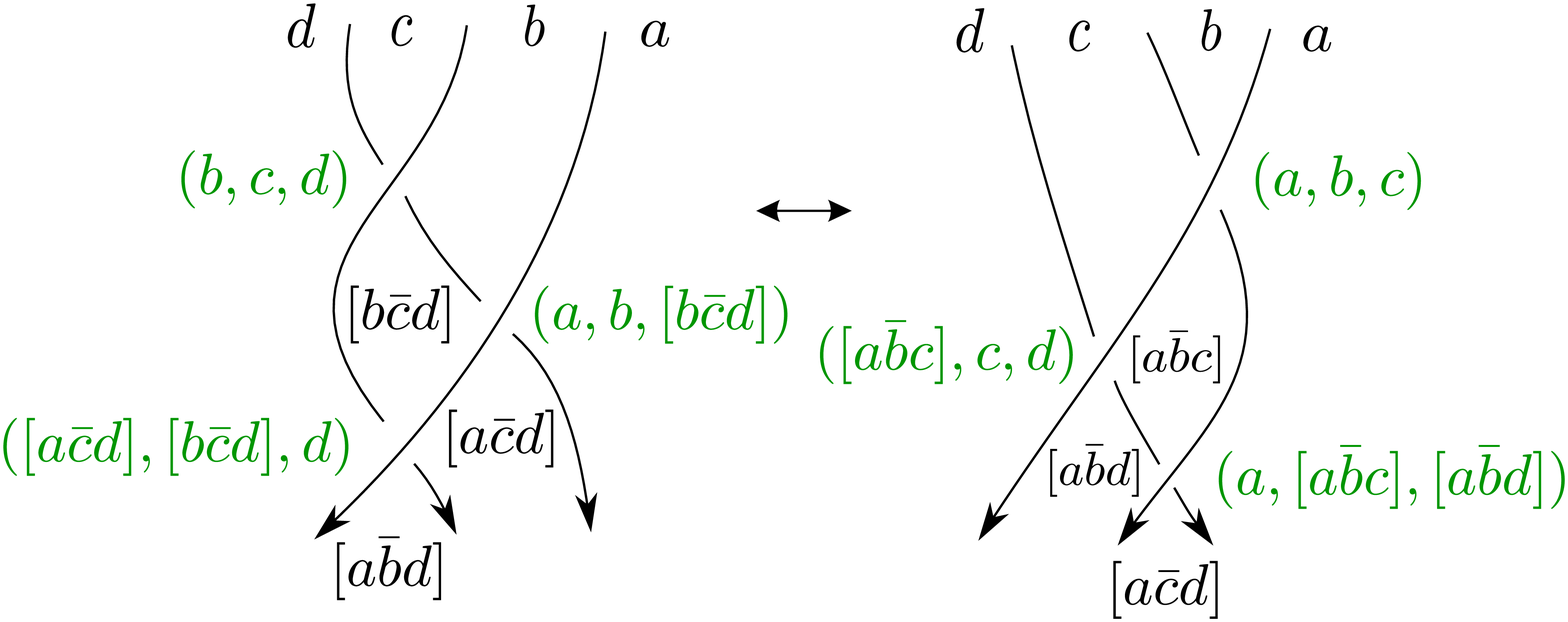}
\caption{}\label{thirddehn}
\end{center}
\end{figure}

The homology class represented by $c_\gamma(D)$ is invariant under all Reidemeister moves. Indeed, the first Reidemeister move contributes (or deletes) a degenerate cycle $(a, b, [b\bar{a}b])$ or $([b\bar{a}b], b, a)$, see Fig. \ref{R1dehn}.
In the second Reidemeister move, because of the sign convention, the contributions cancel, see Fig. \ref{seconddehn}. The invariance under the third Reidemeister move follows from the calculation of the boundary $\partial(a,b,c,d)$ and is depicted in Fig. \ref{thirddehn}.

\begin{figure}
\begin{center}
\includegraphics[height=3.2 cm]{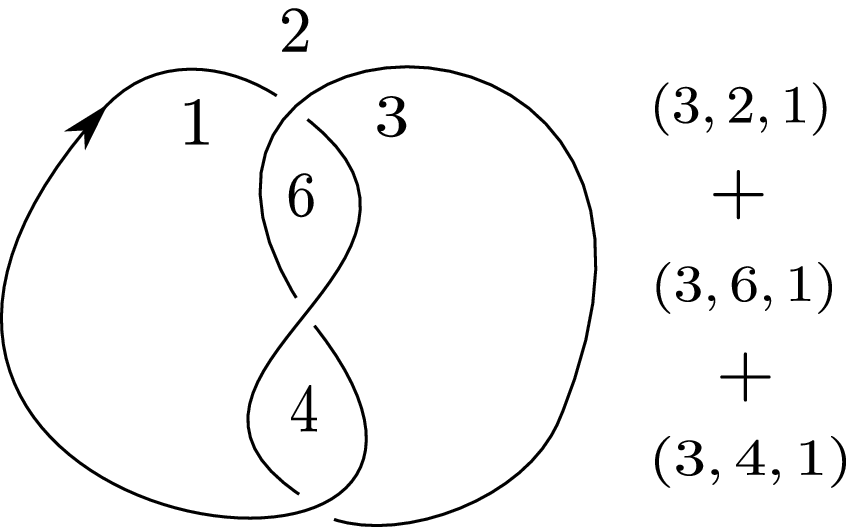}
\caption{}\label{tref223}
\end{center}
\end{figure}

\begin{example}
The trefoil diagram depicted in Fig. \ref{tref223} is colored with elements of ternary group $\mathcal{O}(\triangle(2,2,3))$, whose multiplication tables we included earlier. The cycle assigned to this colored diagram is
\[
(3,2,1)+(3,6,1)+(3,4,1).
\]
It represents $\mathbb{Z}_3$ in the homology of $\mathcal{O}(\triangle(2,2,3))$.
\end{example}

We checked which knots up to 9 crossings have homologically nontrivial colorings
with $\mathcal{O}(\triangle(2,2,3))$ representing torsion elements. All such colorings represented $\mathbb{Z}_3$. We will use the usual knot name from Rolfsen's table and braid notation in which we write, for example [1,-2,3] for $\sigma_1\sigma_2^{-1}\sigma_3$, so the trefoil from Fig. \ref{tref223} is written as [1,1,1]. First we write the name, then the braid notation, after that the number of all colorings, and finally the number of colorings representing $\mathbb{Z}_3$:\\
$3_1$, [ 1, 1, 1 ], 72, 36;\\
$7_4$, [ 1, 1, 2, -1, 2, 2, 3, -2, 3 ], 72, 36;\\
$7_7$, [ 1, -2, 1, -2, 3, -2, 3 ], 72, 36;\\
$8_5$, [ 1, 1, 1, -2, 1, 1, 1, -2 ], 72, 36;\\
$8_{15}$, [ 1, 1, -2, 1, 3, 2, 2, 2, 3 ], 72, 36;\\
$8_{18}$, [ 1, -2, 1, -2, 1, -2, 1, -2 ], 180, 144;\\
$8_{19}$, [ 1, 1, 1, 2, 1, 1, 1, 2 ], 72, 36;\\ 
$8_{21}$, [ 1, 1, 1, 2, -1, -1, 2, 2 ], 72, 36;\\ 
$9_2$, [ 1, 1, 1, 2, -1, 2, 3, -2, 3, 4, -3, 4 ], 72, 36;\\
$9_4$, [ 1, 1, 1, 1, 1, 2, -1, 2, 3, -2, 3 ], 72, 36;\\ 
$9_{10}$, [ 1, 1, 2, -1, 2, 2, 2, 2, 3, -2, 3 ], 72, 36;\\ 
$9_{11}$, [ 1, 1, 1, 1, -2, 1, 3, -2, 3 ], 72, 36;\\ 
$9_{15}$, [ 1, 1, 1, 2, -1, -3, 2, 4, -3, 4 ], 72, 36;\\ 
$9_{16}$, [ 1, 1, 1, 1, 2, 2, -1, 2, 2, 2 ], 72, 36;\\ 
$9_{17}$, [ 1, -2, 1, -2, -2, -2, 3, -2, 3 ], 72, 36;\\ 
$9_{28}$, [ 1, 1, -2, 1, 3, -2, -2, 3, 3 ], 72, 36;\\ 
$9_{29}$, [ 1, -2, -2, 3, -2, 1, -2, 3, -2 ], 72, 36;\\ 
$9_{34}$, [ 1, -2, 1, -2, 3, -2, 1, -2, 3 ], 72, 36;\\ 
$9_{35}$, [ 1, 1, 2, -1, 2, 2, 3, -2, -2, 4, -3, 2, 4, 3 ], 180, 108;\\ 
$9_{37}$, [ 1, 1, -2, 1, 3, -2, -1, -4, 3, -2, 3, -4 ], 180, 72;\\ 
$9_{38}$, [ 1, 1, 2, 2, -3, 2, -1, 2, 3, 3, 2 ], 72, 36;\\ 
$9_{40}$, [ 1, -2, 1, 3, -2, 1, 3, -2, 3 ], 72, 36;\\ 
$9_{46}$, [ 1, -2, 1, -2, 3, 2, -1, 2, 3 ], 180, 72;\\ 
$9_{47}$, [ 1, -2, 1, -2, -3, -2, 1, -2, -3 ], 180, 108;\\ 
$9_{48}$, [ 1, 1, 2, -1, 2, 1, -3, 2, -1, 2, -3 ], 180, 108.
  
One could also use cocycle invariants.
We define the cohomology groups for the ternary group homology in a standard dual way, so the cocycles that are useful for knot invariants are functions
\[
f\colon X\times X\times X\to A,
\]
where $A$ is an abelian group, satisfying
the conditions
\begin{equation*}
\forall_{a,b \in X} \quad f(a, b, [b\bar{a}b])=0,
\end{equation*}
\begin{equation*}
\forall_{a,b \in X} \quad f([b\bar{a}b], b, a)=0,
\end{equation*}
and
\begin{align*}
\forall_{a,b,c,d \in X} \quad 
&f(b,c,d)-f([a\bar{b}c],c,d)-f(a,[a\bar{b}c],[a\bar{b}d])\\
+&f([a\bar{c}d],[b\bar{c}d],d)+f(a,b,[b\bar{c}d])-f(a,b,c)
=0.
\end{align*}

To define cocycle invariants, in the spirit of \cite{CJKLS03}, we consider the cycles corresponding to all possible colorings of regions of a diagram $D$ with a given ternary group $(X, [\ ])$, and evaluate a given cocycle from the cohomology of $(X, [\ ])$ on all these cycles (using the multiplicative notation for the abelian group $A$). The sum of these evaluations is in the group ring $\mathbb{Z}[A]$ and is an invariant of Reidemeister moves.

\section*{Acknowledgements}
I would like to thank Leonid Plakhta and J{\'o}zef H. Przytycki for helpful comments on earlier drafts of this paper.

\bibliography{simplex}
\bibliographystyle{plain}
\end{document}